\documentclass[10pt]{article}

\usepackage[latin1]{inputenc}
\usepackage{amsmath}
\usepackage{amsfonts,amssymb,latexsym,multicol}
\usepackage[francais]{babel}
\usepackage{epsfig,epsf}
\newcounter{fig}
\def\figcaption #1
  {\addtocounter{fig}{1}
   \begin{center}
   {\sc Fig.} \arabic{fig}. {\it #1} 
    \end{center}
    \addtocontents{lof}{\protect \contentsline {figure}{\protect \numberline {\thefig}{\protect \ignorespaces \protect \pit #1 }}{\thepage}}
   }
\input{cyracc.def}

\newtheorem{theo}{Th\'eor\`eme}

\newtheorem{prop}{Proposition}

\newcommand{\implique}{\Rightarrow}

\newcommand{\m}{{\rm m}}

\newcommand{\eps}{\varepsilon}
\newcommand{\fhi}{\varphi}
\newcommand{\ioe}{\leqslant}
\newcommand{\soe}{\geqslant}
\newcommand{\vers}{\rightarrow}
\newcommand{\dist}{{\rm dist}}

\newcommand{\intc}{{\frac{1}{2i\pi}\int \! \!}}
\newcommand{\aintc}{{\frac{1}{2\pi}\int \! \!}}

\newcommand{\demi}{{\frac{1}{2}}}

\newcommand{\Scal}{{\mathcal S}}

\newcommand{\Mgot}{{\mathfrak M}}

\newcommand{\Nat}{{\mathbb N}}

\newcommand{\Real}{{\mathbb R}}
\newcommand{\Com}{{\mathbb C}}

\newcommand{\Res}{{\rm Res}}

\newcommand{\Vect}{{\rm Vect}}

\newcommand{\fin}{\hfill$\Box$}
\newcommand{\dem}{\noindent {\bf D\'emonstration\ }}
\newcommand{\fine}{\tag*{\mbox{$\Box$}}}

\title{Sur un crit\`ere de B\'aez-Duarte pour l'hypoth\`ese de Riemann}
\author{Michel Balazard et Anne de Roton}

\begin{document}
\maketitle

\begin{flushright}
  \textit{Pour Luis B\'aez-Duarte, à l'occasion de son soixante-dixième anniversaire.}
\end{flushright}

\begin{center}
  {\sc Abstract}
\end{center}
\begin{quote}
{\footnotesize Define $e_{n}(t)=\{t/n\} $. Let $d_N$ denote the distance in $L^2(0,\infty ; t^{-2}dt)$ between the indicator function of $[1,\infty[$ and the vector space generated by $e_1, \dots,e_N$. A theorem of B\'aez-Duarte states that the Riemann hypothesis (RH) holds if and only if $d_N \vers 0$ when $N \vers \infty$. Assuming RH, we prove the estimate $$d_N^2 \ioe (\log \log N)^{5/2+o(1)}(\log N)^{-1/2}.$$.}  
\end{quote}

\begin{center}
  {\sc Keywords}
\end{center}
\begin{quote}
{\footnotesize Riemann zeta function, Riemann hypothesis, B\'aez-Duarte criterion, Möbius function. \\MSC classification : 11M26}  
\end{quote}

\section{Position du probl\`eme et \'enonc\'e du r\'esultat principal}

L'étude de la répartition des nombres premiers se ramène à la recherche d'approximations de la fonction
\begin{equation}
  \label{t53}
\chi(x)=[x\soe 1]   
\end{equation}
par des combinaisons linéaires
\begin{equation}
  \label{t54}
  \fhi(x)=\sum_{n=1}^N c_n \{x/n\} \quad (N \in \Nat, \, c_n \in \Real)
\end{equation}
de dilatées de la fonction {\og partie fractionnaire\fg}. Ce fait est connu depuis Tchebychev (cf. \cite{T1852}). En choisissant 
$$
\fhi(x)=-\{x\}+\{x/2\}+\{x/3\}+\{x/5\}-\{x/30\}
$$
il avait observé l'encadrement
$$
\fhi(x)\ioe\chi(x)\ioe\sum_{k \soe 0}\fhi(x/6^k) 
$$
pour en déduire
$$
Ax+O(\log x)\ioe \sum_{n \ioe x} \Lambda(n)\ioe \frac{6}{5}Ax +O(\log^2 x) 
$$
où $\Lambda$ désigne la fonction de von Mangoldt, et
$$
A=\log \frac{2^{1/2}3^{1/3}5^{1/5}}{30^{1/30}}=0,92129202\dots.
$$

On peut préciser la nature de l'approximation de \eqref{t53} par \eqref{t54} équivalente au théorème des nombres premiers
$$
\sum_{n \ioe x} \Lambda(n) \sim  x \quad ( x \vers \infty),
$$
où à l'hypothèse de Riemann
$$
\sum_{n \ioe x} \Lambda(n) =x + O_{\delta}(x^{\demi +\delta}) \quad (x \soe 1, \, \delta >0).
$$
Ainsi, le théorème des nombres premiers est équivalent\footnote{Bien entendu, deux énoncés vrais sont toujours équivalents ; nous renvoyons à \cite{DMC1981} et \cite{BD1993} pour des énoncés précis sur ce sujet.}à l'assertion
$$
\inf_{\fhi} \int_0^{\infty}|\chi(x)-\fhi(x)|\frac{dx}{x^2} =0.
$$
Quant à l'hypothèse de Riemann, B\'aez-Duarte (cf. \cite{BD2003}) a démontré qu'elle équivaut à 
$$
\inf_{\fhi} \int_0^{\infty}|\chi(x)-\fhi(x)|^2\frac{dx}{x^2} =0.
$$
Dans les deux cas, l'infimum est pris sur les $\fhi$ de la forme \eqref{t54}.

\medskip

Nous nous intéressons dans cet article à une forme quantitative de ce critère. Soit $H$ l'espace de Hilbert $L^2(0,\infty ; t^{-2}dt)$ et, pour $\alpha >0$,
$$
e_{\alpha}(t)=\{t/\alpha\} \quad (t>0).
$$

Posons, pour $N$ entier positif,
$$
d_N=\dist_H \bigl (\chi, \Vect(e_1, \dots, e_N)\bigr ).
$$ 

Ainsi, le critère de B\'aez-Duarte affirme que l'hypoth\`ese de Riemann \'equivaut \`a la convergence de $d_N$ vers $0$, quand $N$ tend vers l'infini. 

Examinons maintenant la vitesse de cette convergence. D'une part, Burnol (cf. \cite{B2002}) a d\'emontr\'e que 
$$
d_N^2 \soe \frac{C+o(1)}{\log N}, \quad N \vers +\infty,
$$
o\`u 
$$
C=\sum_{\rho}\frac{\m(\rho)^2}{|\rho|^2},
$$
la somme portant sur les z\'eros non triviaux $\rho$ de la fonction $\zeta$, et $\m(\rho)$ d\'esignant la multiplicit\'e de $\rho$ comme z\'ero de $\zeta$.

Comme
$$
\sum_{\rho}\frac{\m(\rho)}{|\rho|^2}=2+\gamma -\log (4\pi)
$$
(si l'hypothèse de Riemann est vraie, cf. \cite{D2000}, chapter 12, (10), (11)), on en d\'eduit en particulier que
\begin{equation}
  \label{t0}
 d_N^2 \soe \frac{2+\gamma -\log (4\pi)+o(1)}{\log N}, \quad N \vers +\infty. 
\end{equation}

D'autre part, les auteurs de \cite{BDBLS2000} conjecturent l'\'egalit\'e dans \eqref{t0}. Cette conjecture entra\^{\i}ne donc l'hypothèse de Riemann et la simplicit\'e des z\'eros de $\zeta$.

Le comportement asymptotique de $d_N$ est difficile \`a d\'eterminer, m\^eme conditionnellement \`a l'hypothèse de Riemann et d'autres conjectures classiques (simplicit\'e des z\'eros de $\zeta$, conjecture de Montgomery sur la corr\'elation par paires,...). Dans \cite{BD2003}, B\'aez-Duarte donne une d\'emonstration (d\^ue au premier auteur) de la majoration
$$
d_N^2 \ll (\log \log N)^{-2/3}
$$
sous l'hypothèse de Riemann.
Nous am\'eliorons ce r\'esultat dans le pr\'esent travail.
\begin{theo}
  L'hypoth\`ese de Riemann entra\^{\i}ne que
$$
d_N^2 \ll_{\delta} (\log \log N)^{5/2+\delta}(\log N)^{-1/2} \quad (N \soe 3),
$$
pour tout $\delta >0$.
\end{theo}

\medskip

Le plan de notre article est le suivant. Au \S \ref{t18} nous rappelons le rôle de la fonction de M\"obius dans ce problème. Nous y majorons $d_N^2$ par la somme de deux quantités, $I_{N,\eps}$ et $J_{\eps}$, où $\eps$ est un paramètre positif, et nous énonçons les estimations de ces quantités qui permettent de démontrer notre théorème. Le \S \ref{t20} contient une étude de la fonction $\zeta (s)/\zeta (s+\eps)$ nécessaire à la majoration, au \S \ref{t11}, de la quantité $J_{\eps}$. Les \S\S \ref{t55} et \ref{t56} concernent l'estimation des sommes partielles de la série de Dirichlet de l'inverse de la fonction $\zeta$. Cela nous permet de majorer $I_{N,\eps}$ au \S \ref{t19}, concluant ainsi la démonstration.

\medskip

Il apparaîtra clairement que notre travail doit beaucoup à l'article récent \cite{S2008}. Nous remercions son auteur, Kannan Soundararajan, pour une correspondance instructive concernant \cite{S2008}.

\medskip

Le paramètre $\delta$ est fixé une fois pour toutes. On suppose $0<\delta \ioe 1/2$. On pose pour tout nombre complexe $s$
$$
\sigma=\Re s, \quad \tau=\Im s.
$$ 
Les symboles de Bachmann $O$ et de Vinogradov $\ll$ (resp. $\ll_{\delta}$) qui apparaissent sous-entendent toujours des constantes absolues (resp. dépendant uniquement de $\delta$) et effectivement calculables. Enfin nous indiquerons, par les initiales \textit{(HR)} placées au début de l'énoncé d'une proposition, que la démonstration que nous en donnons utilise l'hypothèse de Riemann.

\section{Pertinence de la fonction de M\"obius}\label{t18}

Partant de l'identit\'e
$$
\chi=-\sum_{n \soe 1}\mu(n)e_n
$$
(valable au sens de la convergence simple), B\'aez-Duarte a d'abord montr\'e (cf. \cite{BD1999}) la divergence dans $H$ de la s\'erie du second membre. Il a ensuite propos\'e d'approcher $\chi$ dans $H$ par les sommes
$$
-\sum_{n \ioe N}\mu(n)n^{-\eps}e_n,
$$
où $\eps$ est un paramètre positif, à choisir convenablement en fonction de $N$.

En posant
$$
\nu_{N,\eps}=\Bigl \lVert \chi +\sum_{n \ioe N}\mu(n)n^{-\eps}e_n \Bigr \rVert_H^2,
$$
on a \'evidemment $d_N^2 \ioe \nu_{N,\eps}$ pour $N \soe 1$, $\eps >0$. Posons maintenant pour $N \soe 1$ et $s \in \Com$ :
$$
M_N(s)=\sum_{n \ioe N} \mu(n) n^{-s}.
$$
On sait depuis Littlewood (cf. \cite{L1912}) que l'hypoth\`ese de Riemann entra\^{\i}ne la convergence de $M_N(s)$ vers $\frac{1}{\zeta(s)}$ quand $N$ tend vers l'infini, pour tout $s$ tel que $\Re s >\frac{1}{2}$. Nous allons faire apparaître la différence $M_N-1/\zeta$ pour majorer $\nu_{N,\eps}$.  
\begin{prop}\label{t24}
  Pour $N \soe 1$ et $\eps >0$, on a
$$
\nu_{N,\eps} \ioe 2I_{N,\eps} +2J_{\eps},
$$
o\`u
$$
I_{N,\eps} = \frac{1}{2\pi} \int_{\sigma =1/2} |\zeta (s)|^2 |M_N(s+\eps)-\zeta(s+\eps)^{-1}|^2\frac{d\tau}{|s|^2} \quad \text{ et } \quad J_{\eps} = \frac{1}{2\pi} \int_{\sigma =1/2} \left \vert \frac{\zeta (s)}{\zeta (s+\eps)}-1 \right \vert ^2 \frac{d\tau}{|s|^2}.
$$
\end{prop}
\dem

La transformation de Mellin associe \`a toute fonction $f \in H$ une fonction $\Mgot{f}$, d\'efinie pour presque tout $s$ tel que $\sigma=1/2$ par la formule
$$
\Mgot{f}(s)=\int_0^{+\infty} f(t)t^{-s-1}dt
$$
(o\`u $\int_0^{+\infty}$ signifie $\lim_{T \vers +\infty}\int_{1/T}^T$).

De plus, le th\'eor\`eme de Plancherel affirme que $f \mapsto \Mgot{f}$ est un op\'erateur unitaire entre $H$ et \\$L^2(\frac{1}{2}+i\Real , d\tau /2\pi)$, espace que nous noterons simplement $L^2$. Comme
$$
\Mgot{e_{\alpha}}(s)=\alpha^{-s}\frac{\zeta (s)}{-s}, \quad \quad \Mgot{\chi}(s)=\frac{1}{s},$$
on a
\begin{align*}
\nu_{N,\eps}&=\Bigl \lVert \chi +\sum_{n \ioe N}\mu(n)n^{-\eps}e_n \Bigr \rVert_H^2\\
&= \Bigl \lVert \frac{1}{s} +\sum_{n \ioe N}\mu(n)n^{-\eps} n^{-s}\frac{\zeta (s)}{-s} \Bigr \rVert_{L^2}^2\\
&= \frac{1}{2\pi}  \int_{\sigma =1/2} |1-\zeta (s)M_N(s+\eps)|^2\frac{d\tau}{|s|^2}\\
& \ioe \frac{1}{\pi} \int_{\sigma =1/2} \left \vert 1-\frac{\zeta (s)}{\zeta (s+\eps)} \right \vert ^2 \frac{d\tau}{|s|^2} +\frac{1}{\pi} \int_{\sigma =1/2} \left \vert \frac{\zeta (s)}{\zeta (s+\eps)}- \zeta (s)M_N(s+\eps) \right \vert ^2 \frac{d\tau}{|s|^2}\\
& \text{\footnotesize (o\`u l'on a utilis\'e l'in\'egalit\'e $|a+b|^2 \ioe 2 (|a|^2+|b|^2)$)}\\
&=2J_{\eps} + 2I_{N,\eps}.\fine
\end{align*}

Observons que la proposition \ref{t24} ne dépend pas de l'hypothèse de Riemann, mais que les quantités $I_{N,\eps}$ et $J_{\eps}$ pourraient être infinies si elle était fausse.

Dans \cite{BD2003}, B\'aez-Duarte démontre (sous l'hypothèse de Riemann) que $I_{N,\eps}$ tend vers $0$ quand $N$ tend vers l'infini (pour tout $\eps >0$ fixé), et que $J(\eps)$ tend vers $0$ quand $\eps$ tend vers $0$. On a donc bien $d_N=o(1)$.

La version quantitative donnée dans \cite{BD2003} repose sur les estimations
$$
J(\eps) \ll \eps^{2/3} \quad (0 <\eps \ioe 1/2),
$$ 
et
$$
I_{N,\eps} \ll N^{-2\eps/3} \quad (c/\log \log N \ioe \eps \ioe 1/2),
$$
où $c$ est une constante positive absolue. Nous démontrons ici les deux propositions suivantes.
\begin{prop}\label{t27}(HR)
 On a $J_{\eps} \ll \eps$.
\end{prop}

\begin{prop}\label{t57}(HR)
Soit $\delta >0$. Pour $N\soe N_0(\delta)$  et $\eps \soe 25(\log \log N)^{5/2+\delta}(\log N)^{-1/2}$, on a
$$
I_{N,\eps} \ll N^{-\eps/2}.
$$  
\end{prop}

Le choix $\eps =25 (\log \log N)^{5/2+\delta}(\log N)^{-1/2}$ donne le théorème.

\section{\'Etude du quotient $\zeta (s)/\zeta (s+\eps)$}\label{t20}

Dans ce paragraphe, nous \'etudions, sous l'hypoth\`ese de Riemann, le comportement de la fonction $\zeta (s)/\zeta (s+\eps)$ dans le demi-plan $\sigma \soe 1/2$, quand $\eps$ tend vers $0$. Afin de pr\'eciser, sur certains points, l'expos\'e de Burnol dans \cite{B2003}, nous utilisons le produit de Hadamard de $\zeta(s)$ et majorons chaque facteur de $\zeta (s)/\zeta (s+\eps)$.

Nous supposons $0<\eps \ioe 1/2$.

\begin{prop}\label{t1}(HR)
 On a les estimations suivantes.
 \begin{align*}
   (i) \quad\quad \left \lvert \frac{\zeta (s)}{\zeta (s +\eps)}\right \rvert^2 &\ll |s|^{\eps} \quad (\sigma=1/2) ;\\ 
(ii) \quad \quad \left \lvert \frac{\zeta (s)}{\zeta (s +\eps)}\right \rvert^2 &\ioe 1+O(\eps|s|^{1/2})\quad (\sigma=1/2);\\
(iii) \quad \frac{\zeta (s)/\zeta (s +\eps)}{s(1-s)} &\ll \frac{|s|^{\eps/2}}{|s-1|^2} \quad(\sigma \soe 1/2, \, s \not =1).
 \end{align*}
\end{prop}
\dem

Si l'on pose 
$$
\xi (s) = \frac{1}{2} s(s-1)\pi^{-s/2}\Gamma (s/2) \zeta (s),
$$
on a
$$
\xi (s)=\prod_{\rho} \left ( 1-\frac{s}{\rho}\right ),
$$
o\`u le produit porte sur les z\'eros non triviaux $\rho$ de la fonction $\zeta$, et doit \^etre calcul\'e par la formule $\prod_{\rho}=\lim_{T \vers +\infty} \prod_{|\gamma | \ioe T}$ (on pose $\rho = \beta + i\gamma$). Par conséquent
 \begin{equation}
   \label{t7}
 \frac{\zeta (s)}{\zeta (s +\eps)}=\pi^{-\eps/2}\frac{(s+\eps)(s+\eps-1)}{s(s-1)}\frac{\Gamma \bigl ( (s+\eps)/2 \bigr )}{\Gamma (s/2)}\prod_{\rho} \frac{s-\rho}{s+\eps -\rho}.  
 \end{equation} 

\medskip

Examinons successivement les facteurs apparaissant dans \eqref{t7}. On a d'abord $\pi^{-\eps/2} <1$. Ensuite, on a
\begin{align}
\left \lvert\frac{(s+\eps)(s+\eps-1)}{s(s-1)} \right \rvert &\ll   \left \lvert \frac{s}{s-1} \right \rvert\quad (\sigma \soe 1/2, \, s \not =1),\label{t13}\\
 \left \lvert \frac{(s+\eps)(s+\eps-1)}{s(s-1)}  \right \rvert &\ioe \exp\bigl(O(\eps/|s|)\bigr) \quad (\sigma = 1/2).\label{t59}
\end{align}

Pour le quotient des fonctions $\Gamma$ apparaissant dans la formule \eqref{t7}, on dispose de l'in\'egalit\'e suivante, qui résulte de la formule de Stirling complexe.

\begin{equation}  \label{t8}
\left \lvert \frac{\Gamma \bigl ( (s+\eps)/2 \bigr )}{\Gamma (s/2)}\right \rvert \ioe |s/2|^{\eps/2}  \exp\bigl(O(\eps/|s|)\bigr ) \quad (\sigma \soe 1/2).
\end{equation}

Pour majorer le produit infini apparaissant dans \eqref{t7}, on utilise l'in\'egalit\'e 
$$
\Bigl \lvert  \frac{s-\rho}{s+\eps -\rho} \Bigr \rvert <1, \quad \sigma \soe \beta, \quad \eps >0,
$$
qui donne par conséquent (sous l'hypothèse de Riemann)
\begin{equation} 
  \label{t10}
 \left \lvert \prod_{\rho} \frac{s-\rho}{s+\eps -\rho} \right \rvert <1  \quad (\sigma \soe 1/2).
\end{equation}

Notons ensuite les inégalités 
\begin{equation}
  \label{t15}
 \exp \bigl (\eps \log x/2+O(\eps/x)\bigr )\ll(x/2)^{\eps}, 
\end{equation}
et
\begin{equation}
  \label{t58}
\exp \bigl (\eps \log x/2+O(\eps/x)\bigr ) \ioe 1 +O( \eps x^{1/2}), 
\end{equation}
valables pour $x\soe 1/2$.

L'estimation \textit{(i)} résulte alors de \eqref{t7}, \eqref{t59}, \eqref{t8}, \eqref{t10} et \eqref{t15} ; l'estimation \textit{(ii)} de \eqref{t7}, \eqref{t59}, \eqref{t8}, \eqref{t10} et \eqref{t58}, et l'estimation \textit{(iii)} de \eqref{t7}, \eqref{t13}, \eqref{t8}, \eqref{t10} et \eqref{t15}.\fin

\section{Majoration de $J_{\eps}$}\label{t11}

On suppose, comme au \S \ref{t20}, que $\eps$ v\'erifie $0<\eps \ioe 1/2$. On pose
$$
K_{\eps}=\frac{1}{2\pi} \int_{\sigma =1/2} \left \vert \frac{\zeta (s)}{\zeta (s+\eps)} \right \vert ^2 \frac{d\tau}{|s|^2} \quad \text{et} \quad L_{\eps} =\frac{1}{2\pi} \int_{\sigma =1/2} \frac{\zeta (s)}{\zeta (s+\eps)} \frac{d\tau}{|s|^2}, 
$$
de sorte que
\begin{equation}\label{t60}
J_{\eps}=K_{\eps}-2L_{\eps}+1.  
\end{equation}

Pour majorer $J_{\eps}$, nous allons calculer exactement $L_{\eps}$ à l'aide du théorème des résidus, et majorer $K_{\eps}$ en utilisant les résultats du paragraphe précédent.

\begin{prop}\label{t3}(HR)
  On a 
  \begin{align*}
L_{\eps}&=\frac{\gamma -1}{\zeta(1+\eps)} -\frac{\zeta ' (1+\eps)}{\zeta^2(1+\eps)} \\
&=1-(\gamma +1) \eps +O(\eps^2).    
  \end{align*}
\end{prop}
\dem

On a 
\begin{align*}
L_{\eps} &=  \frac{1}{2\pi} \int_{\sigma =1/2} \frac{\zeta (s)}{\zeta (s+\eps)} \frac{d\tau}{|s|^2}\\
&= \frac{1}{2\pi i} \int_{\sigma =1/2}Q(s)ds,
\end{align*}
où
$$
Q(s)=\frac{\zeta (s)/\zeta (s +\eps)}{s(1-s)}.
$$

Soit $\Pi$ le demi-plan $\sigma \soe \frac{1}{2}$, et $\Delta$ la droite $\sigma = \frac{1}{2}$. La fonction $Q$ est m\'eromorphe dans $\Pi$, holomorphe sur $\Delta$. Dans $\Pi$ elle a un unique p\^ole, double, en $s=1$ o\`u son r\'esidu vaut
$$
\frac{1-\gamma }{\zeta(1+\eps)} +\frac{\zeta'(1+\eps)}{\zeta^2(1+\eps)}.
$$
D'après la proposition \ref{t1}, \textit{(iii)}, on a $sQ(s) \vers 0$ uniform\'ement quand $|s| \vers +\infty$, $s \in \Pi$, et
$$
\int_{\Delta} |Q(s)| \cdot |ds| < +\infty.
$$

Nous sommes donc en situation d'appliquer une proposition classique du calcul des résidus (cf. par exemple \cite{WW1927}\S 6.22) pour en déduire
\begin{align*}
  L_{\eps} &=-\Res \left ( \frac{\zeta (s)}{\zeta (s+\eps)}\cdot \frac{1}{s(1-s)} \right )\Big \vert _{s=1}\\
&= \frac{\gamma -1}{\zeta(1+\eps)} -\frac{\zeta ' (1+\eps)}{\zeta^2(1+\eps)}. 
\end{align*}

Cette dernière quantité vaut
$$
1-(\gamma +1) \eps +O(\eps^2)
$$
puisque
\begin{equation*}
\frac{1}{\zeta (1+\eps)} = \eps -\gamma \eps^2+O(\eps^3).\fine  
\end{equation*}
  
\medskip

Nous sommes maintenant en mesure de démontrer l'estimation $J_{\eps} \ll \eps$, objet de la proposition \ref{t27}. En intégrant l'inégalité \textit{(ii)} de la proposition \ref{t1} sur la droite $\sigma=1/2$ avec la mesure $d\tau/|s^2|$, on obtient 
$$
K_{\eps}-1 \ll \eps.
$$
Le résultat découle alors de \eqref{t60} et de la proposition \ref{t3}.

\medskip
 
En considérant la contribution à $J_{\eps}$ d'un voisinage de l'ordonnée d'un zéro simple de $\zeta$ (par exemple $\gamma_1=14,1347\dots$), on peut montrer inconditionnellement que $J_{\eps} \gg \eps$. Il serait intéressant de préciser le comportement asymptotique de $J_{\eps}$ quand $\eps$ tend vers $0$.

\section{Quelques propriétés de la fonction $\zeta$ sous l'hypothèse de Riemann}\label{t55}

Afin d'établir la majoration de la proposition \ref{t57}, nous allons étudier $M_N(s+\eps)$. Pour cela, nous allons utiliser la méthode inventée par Maier et Montgomery dans l'article \cite{MM2008}, dévolu à $M_N(0)=M(N)$. Ils y démontrent que
$$
M(N)=\sum_{n \ioe N} \mu(n) \ll \sqrt{N}\exp \bigl ( (\log N)^{39/61}\bigr )
$$
sous l'hypothèse de Riemann. Leur approche a été ensuite perfectionnée par Soundararajan (cf. \cite{S2008}), qui a obtenu l'estimation
$$
M(N) \ll \sqrt{N}\exp \bigl ( (\log N)^{1/2} (\log \log N)^{14}\bigr ),
$$
toujours sous l'hypothèse de Riemann. La méthode de Soundararajan donne en fait
$$
M(N) \ll_{\delta} \sqrt{N}\exp \bigl ( (\log N)^{1/2} (\log \log N)^{5/2+\delta}\bigr ),
$$
pour tout $\delta$ tel que $0<\delta \ioe 1/2$. Nous allons maintenant rappeler les éléments de la méthode de Soundararajan qui seront utilisés dans notre argumentation, avec les quelques modifications qui permettent d'obtenir l'exposant $5/2+\delta$. On trouvera les démonstrations dans l'article \cite{S2008} (cf. aussi \cite{BR2008} pour un exposé détaillé des modifications).

\subsection{Ordonnées $V$-typiques}\label{t36}

L'évaluation de $M_N(s+\eps)$ grâce à la formule de Perron fera appel à un contour sur lequel les grandes valeurs de $|\zeta (z)|^{-1}$ seront aussi rares que possible. Pour quantifier cette rareté, Soundararajan a introduit la notion suivante.

Soit $T$ assez grand\footnote{Ici et dans la suite, cela signifie que $T \soe T_0(\delta)$, quantité effectivement calculable, et dépendant au plus de $\delta$.} et $V$ tel que $(\log \log T)^2 \ioe V \ioe \log T/\log\log T$. Un nombre réel $t$ est appelé une \textbf{ordonnée $V$-typique de taille $T$} si

$\bullet$ $T \ioe t \ioe 2T$ ;

\textit{(i)} pour tout $\sigma \soe 1/2$, on a
$$
\Bigl \lvert \sum_{n \ioe x} \frac{\Lambda (n)}{n^{\sigma +it}\log n}\frac{\log (x/n)}{\log x} \Bigr \rvert \ioe 2V, \quad \text{où $x=T^{1/V}$} ;
$$

\textit{(ii)} tout sous-intervalle de $[t-1,t+1]$ de longueur $2\pi\delta V/\log T$ contient au plus $(1+\delta)V$ ordonnées de zéros de $\zeta$ ;

\textit{(iii)} tout sous-intervalle de $[t-1,t+1]$ de longueur $2\pi V/\bigl ((\log V)\log T\bigr )$ contient au plus $V$ ordonnées de zéros de $\zeta$.

Si $t\in [T,2T]$ ne vérifie pas l'une des assertions \textit{(i), (ii), (iii)}, on dira que $t$ est une \textbf{ordonnée $V$-atypique de taille $T$}.

L'apport de cette définition à l'estimation de $M_N(s+\eps)$ via la formule de Perron (\S \ref{t56} ci-dessous) est contenu dans l'énoncé suivant (proposition 9 de \cite{BR2008}).

\begin{prop}(HR)\label{t43}
Soit $t$ assez grand, et $x\soe t$. Soit $V'$ tel que $(\log \log t)^2 \ioe V' \ioe (\log t/2) /(\log \log t/2)$. On suppose que $t$ est une ordonnée $V'$-typique (de taille $T'$). Soit $V \soe V'$.

Alors 
$$
|x^z\zeta(z)^{-1}| \ioe \sqrt{x} \exp \bigl ( V\log (\log x/\log t) +(2+3\delta)V\log \log V\bigr) \quad\quad (V' \ioe (\Re z -1/2)\log x \ioe V, \quad |\Im z|=t).
$$
\end{prop}

\subsection{Majoration de l'écart entre le nombre de zéros de la fonction $\zeta$ et sa moyenne, dans un intervalle de la droite critique}\label{t29} 

La proposition suivante (cf. \cite{BR2008}, proposition 15) donne une majoration de l'écart entre le nombre d'ordonnées de zéros de $\zeta$ dans l'intervalle $]t-h,t+h]$ et sa valeur moyenne $(h/\pi)\log(t/2\pi)$. Cet encadrement est exprimé au moyen d'un paramètre $\Delta$, et met notamment en jeu un polynôme de Dirichlet de longueur $\exp 2\pi \Delta$.
\begin{prop}\label{t34}(HR)
Soit $\Delta \soe 2$ et $h>0$. Il existe des nombres réels $a(p)=a(p,\Delta,h)$ ($p$ premier, $p \ioe e^{2\pi\Delta}$) vérifiant 

$\bullet$ $|a(p)| \ioe 4$ pour $p \ioe e^{2\pi \Delta}$ ;

$\bullet$ pour tout $t$ tel que $t \soe \max(4,h^2)$, on a
$$
N(t+h)-N(t-h)- 2h\frac{\log t/2\pi}{2\pi}\ioe  \frac {\log t}{2\pi\Delta}+ \sum_{p\ioe e^{2\pi\Delta}} \frac{a(p)\cos(t\log p)}{p^{\demi}}+O (\log  \Delta ).
$$  
\end{prop}

\medskip

Lorsqu'on majore trivialement le polynôme de Dirichlet qui intervient dans cette proposition, on obtient le résultat suivant, dû à Goldston et Gonek (cf. \cite{GG2007}). Notre énoncé est légèrement plus précis que celui de \cite{GG2007}.
\begin{prop}\label{t42}
 Soit $t$ assez grand et $0<h \ioe \sqrt{t}$. On a
$$
N(t+h)-N(t-h)- (h/\pi)\log (t/2\pi)\ioe (\log t)/2\log \log t + \bigl (1/2 +o(1)\bigr ) \log t\log \log \log t/(\log \log t)^2.
$$ 
\end{prop}
\dem

On a
\begin{align*}
  \Bigl \lvert  \sum_{p\ioe e^{2\pi\Delta}} \frac{a(p)\cos t\log p}{p^{\demi}} \Bigr \rvert& \ll \sum_{p\ioe e^{2\pi\Delta}} \frac{1}{\sqrt{p}}\\
&\ll \frac{ e^{\pi\Delta}}{\Delta}.
\end{align*}
On choisit $\Delta = \frac{1}{\pi}\log(\log t/\log \log t)$ et on vérifie alors que 
\begin{equation*}
 \frac {\log t}{2\pi\Delta}+O(e^{\pi\Delta}/\Delta) +O (\log  \Delta )= (\log t)/2\log \log t +  \bigl (1/2 +o(1)\bigr ) \log t\log \log \log t/(\log \log t)^2.\fine 
\end{equation*}

\medskip

La proposition suivante est une variante un peu plus précise de la première assertion de la Proposition 4 de \cite{S2008}.
\begin{prop}\label{t41}
  Soit $T$ assez grand, et $V$ tel que
$$
\demi + \Bigl (\demi +\delta\Bigr ) \log \log \log T/\log \log T \ioe V \log \log T/ \log T \ioe 1.
$$
Alors toute ordonnée $t \in [T,2T]$ est $V$-typique.
\end{prop}
\dem

Il faut vérifier les critères \textit{(i), (ii), (iii)} de la définition d'une ordonnée $V$-typique.

Pour \textit{(i)}, on a  pour $\sigma \soe 1/2$, $t \in \Real$, et $x=T^{1/V}$,
\begin{align*}
\Bigl \lvert \sum_{n \ioe x} \frac{\Lambda (n)}{n^{\sigma +it}\log n}\frac{\log (x/n)}{\log x} \Bigr \rvert & \ioe  \sum_{n\ioe x}\frac{\Lambda (n)}{\sqrt{n}\log n}\frac{\log (x/n)}{\log x}\\
& \ll  \frac{\sqrt{x}}{(\log x)^2}\\
& \ll  \frac{\log T}{(\log \log T)^2} \quad \text{\footnotesize (car $x=T^{1/V} \ioe (\log T)^2$)}\\
& = o(V). 
\end{align*}

Pour \textit{(ii)} on a, avec $t' \in[t-1,t+1]$ et $h=\pi\delta V/\log T$ :
\begin{align*}
N(t'+h)-N(t'-h) &\ioe  (h/\pi)\log (t'/2\pi) +\demi \log t'/\log \log t' + \bigl (1/2 +o(1)\bigr ) \log t'\log \log \log t'/(\log \log t')^2\\
 & \quad \text{\footnotesize (proposition \ref{t42})}\\ 
& \ioe  (h/\pi)\log T + \demi \log T/\log \log T + (1/2+\delta)\log T\log \log \log T/(\log \log T)^2\\
& \ioe (1+\delta)V.
\end{align*}

Pour \textit{(iii)} on a, avec $t' \in[t-1,t+1]$ et $h=\pi V/\bigl ((\log V)\log T\bigr )$ :
\begin{align*}
N(t'+h)-N(t'-h) &\ioe  (h/\pi)\log (t'/2\pi) +\demi \log t'/\log \log t' + \bigl (1/2 +o(1)\bigr ) \log t'\log \log \log t'/(\log \log t')^2\\
 & \ioe  \frac{V}{\log V} + \demi \log T/\log \log T +  \bigl (1/2 +o(1)\bigr )\log T\log \log \log T/(\log \log T)^2\\
& \ioe \demi \log T/\log \log T + (1/2+\delta)\log T\log \log \log T/(\log \log T)^2\\
& \ioe V.\fine
\end{align*}

\section{Approximation de l'inverse de la fonction $\zeta$ par ses sommes partielles}\label{t56}

Le but de ce paragraphe est la démonstration de la proposition suivante.
\begin{prop}\label{t52}
  Soit $N$ assez grand et $\eps \soe 25(\log \log N)^{5/2+6\delta}(\log N)^{-1/2}$. Alors, pour $|\tau| \ioe N^{3/4}$,on a
$$
\zeta (s+\eps)^{-1}-M_N(s+\eps) \ll N^{-\eps/4}(1+|\tau|)^{1/2-\beta(\tau)},
$$
où $\beta (\tau)=\frac{\log\log\log(16+|\tau|)}{2\log\log(16+|\tau|)}$. 
\end{prop}

Elle résultera de diverses estimations, valables uniformément quand $\tau$ et $\eps$ appartiennent à certains intervalles définis en termes de $N$, longueur du polynôme de Dirichlet $M_N$, approximant la fonction $\zeta^{-1}$. Pour plus de clarté dans l'exposé, nous développons séparément les analyses relatives aux deux paramètres $\tau$ et $\eps$. Nous commençons par l'étude de  
$$
M_N(i\tau) =\sum_{n \ioe N}\mu(n)n^{-i\tau},
$$
pour $\tau \in \Real$.

\subsection{Estimation de $M_N(i\tau)$ pour les petites valeurs de $|\tau|$}

Commençons par le résultat obtenu par sommation partielle à partir de la majoration de Soundararajan (cf. \cite{S2008} et \cite{BR2008})
$$
M(x) =\sum_{n \ioe x}\mu (n) \ll \sqrt{x}\exp C(\log x), \quad x \soe 3,
$$
où $C(u)=u^{1/2}(\log u)^{5/2+\delta}$. Observons que $C'(u)=O(1)$, $u \soe 1$.
\begin{prop}\label{t48}
  On a uniformément
$$
M_N(i\tau) \ll (1+|\tau|)  \sqrt{N}\exp C(\log N), \quad N \soe 3, \quad \tau \in \Real.
$$
\end{prop}

 La démonstration (standard) est laissée au lecteur. Pour aller plus loin, nous allons appliquer la formule de Perron et suivre la démarche de Soundararajan dans \cite{S2008}.

\subsection{Estimation de $M_N(i\tau)$ pour les grandes valeurs de $|\tau|$}

Nous utiliserons la majoration simple suivante.
\begin{prop}\label{t44}
Pour $0<\delta \ioe 1/12$, $N$ assez grand et
$$
\exp\bigl (3(\log N)^{1/2}(\log \log N)^{5/2+6\delta}\bigr ) \ioe |\tau| \ioe N^{3/4},
$$
on a
$$
M_N(i\tau) \ll N^{1/2}|\tau|^{1/2-\kappa(\tau)}, 
$$
où $\kappa (\tau)=\demi\log\log\log |\tau|/\log\log |\tau|$.
\end{prop}

\dem
Dans toute la démonstration, $N$ sera supposé assez grand.

\subsubsection*{Première étape : formule de Perron}

La première étape de la démonstration consiste à appliquer la formule de Perron à la hauteur 
$N_1= 2^{\lfloor \log N/\log 2\rfloor}$ (le choix d'une puissance de $2$ simplifie l'exposé de \cite{BR2008}), ce qui pour $\tau \in \Real$ donne
\begin{align*}
M_N(i\tau) &=\frac{1}{2\pi i }\int_{1+1/\log N-iN_1}^{1+1/\log N+iN1}\zeta(z+i\tau)^{-1}\frac{N^z}{z}dz + O(N\log N_1/N_1)\\
&= \frac{1}{2\pi i }\int_{1+1/\log N-i(N_1-\tau)}^{1+1/\log N+i(N_1+\tau)}\zeta(z)^{-1}\frac{N^{z-i\tau}}{z-i\tau}dz+O(\log{N})\\  
\end{align*}
Supposons maintenant que $|\tau| \ioe N/5$ et remplaçons l'intégrale par $N^{-i\tau}B_N$, où
$$
B_N =B_N(i\tau)=\frac{1}{2\pi i }\int_{1+1/\log N-iN_1}^{1+1/\log N+iN_1}\zeta(z)^{-1}\frac{N^{z}}{z-i\tau}dz.
$$
L'erreur commise est alors majorée par
$$\frac{1}{2\pi}\int_{N_1-|\tau| \ioe |\Im z| \ioe N_1+|\tau|}|\zeta(z)^{-1}| \Bigl \lvert\frac{N^{z}dz}{z-i\tau}\Bigr \rvert \quad (\Re{z}=1+1/\log N).$$
Or $|\zeta(z)^{-1}| \ll \log N$ si $\Re z=1+1/\log N$ et $|z-i\tau| \gg N$ si $N_1-|\tau| \ioe |\Im z| \ioe N_1+|\tau|$, donc l'erreur est $O(|\tau| \log N)$.\\
Pour $N \soe 3$ et $|\tau|\ioe N/5$ on a donc montré
\begin{equation}\label{t45}
M_N(i\tau)=N^{-i\tau}B_N+O\bigl ((1+|\tau|) \log N \bigr ).
\end{equation}

\subsubsection*{Deuxième étape : déformation du chemin d'intégration}

Pour majorer $|B_N|$, nous allons remplacer le segment d'intégration $[1+1/\log N-iN_1,1+1/\log N +iN_1]$ par  une variante $\Scal_N$ du chemin défini par Soundararajan dans \cite{S2008}, chemin sur lequel les grandes valeurs de l'intégrande sont rares. Nous commençons par une description de $\Scal_N$. Nous posons
$$
\kappa=\lfloor  (\log N)^{1/2}(\log \log N)^{5/2} \rfloor, \quad K =\lfloor \log N/\log 2\rfloor.
$$
Nous posons également $T_k=2^{k}$ pour $\kappa \ioe k \ioe K$, et $N_0=T_{\kappa}$ (on a $N_1=T_K$).

\medskip

Le chemin $\Scal_N$ est symétrique par rapport à l'axe réel, et constitué de segments verticaux et horizontaux. Nous décrivons seulement la partie de $\Scal_N$ située dans le demi-plan $\Im z \soe 0$.

$\bullet$ Il y a d'abord un segment vertical $[1/2+1/\log N,1/2+1/\log N +iN_0]$.

$\bullet$ Pour chaque $k$ tel que $\kappa \ioe k < K$, on considère les entiers $n$ de l'intervalle $[T_k,2T_k[$. On définit alors $V_n$ comme le plus petit entier de l'intervalle $[(\log \log T_k)^2,\log T_k/\log \log T_k]$ tel que tous les points de $[n,n+1]$ soient $V_n$-typiques de taille $T_k$. L'existence de $V_n$ est garantie par la proposition \ref{t41}. On a même 
$$
V_n \ioe \demi  \log n/\log \log n +(1/2+\delta)\log n (\log \log \log n )/(\log \log n)^2 +1.
$$ 
On inclut alors dans $\Scal_N$ le segment vertical $[1/2+V_n/\log N +in,1/2+V_n/\log N +i(n+1)]$

Il y a enfin des segments horizontaux reliant tous ces segments verticaux :

$\bullet$ le segment $[1/2+1/\log N +iN_0,1/2+V_{N_0}/\log N +iN_0]$ ;

$\bullet$ les segments $[1/2+V_n/\log N +i(n+1),1/2+V_{n+1}/\log N +i(n+1)]$, $N_0 \ioe n \ioe T_{K}-2$ ;

$\bullet$ le segment $[1/2+V_{N_1-1}/\log N +iN_1,1+1/\log N +iN_1]$.

\medskip

D'après le théorème de Cauchy, on a
$$
B_N=\intc_{\Scal_N}\zeta(z)^{-1}\frac{N^{z}}{z-i\tau}dz.
$$

\subsubsection*{Troisième étape : \'evaluation de $B_N$}
Lorsque $|z-i\tau|$ n'est pas trop petit devant $|z|$, nous pouvons utiliser les estimations de \cite{S2008} et \cite{BR2008}. Nous définissons donc $\Scal_{N,\tau}$ comme la partie de $\Scal_N$ où $  |(\Im z -\tau)/\tau| \ioe 1/4$ ($\tau \not =0$).

 Si $z \in \Scal_{N} \setminus \Scal_{N,\tau}$, on a $|z-i\tau|\gg |z|$. Par conséquent (cf. \cite{S2008} et \cite{BR2008}), pour $N \soe 3$ et $\tau \in \Real$, on a
\begin{align}\label{t46}
  \Bigl \lvert B_N-\intc_{\Scal_{N,\tau}}\zeta(z)^{-1}\frac{N^{z}}{z-i\tau}dz \Bigr \rvert &\ll\int_{ \Scal_{N}} \Bigl \lvert\frac{\zeta(z)^{-1}N^{z}dz}{z}\Bigr \rvert\nonumber\\
&\ll \sqrt{N} \exp\bigl ((\log N)^{1/2}(\log \log N)^{5/2+6\delta}\bigr ). 
\end{align}

\medskip

Il nous reste à majorer la contribution de $\Scal_{N,\tau}$.

Supposons $ \sqrt{2}N_0 \ioe |\tau|  \ioe \frac{1}{\sqrt{2}}N_1$. Par symétrie, on peut également supposer $\tau >0$. On a
$$
\left\vert\intc_{\Scal_{N,\tau}}\zeta(z)^{-1}\frac{N^{z}}{z-i\tau}dz\right\vert \ioe \sup_{z \in \Scal_{N,\tau}}|\zeta(z)^{-1}N^{z}| \left(\aintc_{\Scal_{N,\tau}} \Bigl \lvert\frac{dz}{z-i\tau}\Bigr \rvert\right).
$$
Observons que si $z \in \Scal_N$ et $\Im z  \soe N_0$, alors $z$ se trouve sur un des segments horizontaux et verticaux décrits ci-dessus. Sur les deux segments (horizontal et vertical) de $\Scal_{N,\tau}$ situés dans la bande $n< \Im z \ioe n+1$, on a $|z-i\tau|^{-1} \ll (1+|n-\tau|)^{-1}$, donc l'intégrale est en $O( \log \tau)$.

\medskip

Pour majorer $|\zeta(z)^{-1}N^{z}|$, nous utilisons la proposition \ref{t43}.
En posant $n=\lceil \Im z \rceil -1$, on peut écrire
$$
V' \ioe (\Re z -1/2)\log N \ioe V,
$$
avec $(V,V')=(V_n,V_n)$ dans le cas vertical et $(V_{n+1},V_n)$ ou $(V_n,V_{n+1})$ dans le cas horizontal ($\Im{z}=n+1$), et $\Im z$ $V'$-typique (de taille correspondante). On peut donc bien appliquer la proposition \ref{t43} pour obtenir 
$$
|\zeta(z)^{-1}N^{z}| \ioe \sqrt{N} \exp \bigl ( V\log (\log N/\log \Im z) +(2+3\delta)V\log\log V \bigr ).
$$

Maintenant, si $z \in \Scal_{N,\tau}$, on a
$$
\tau \sqrt{2} \soe \Im z \soe \tau/ \sqrt{2} \soe N_0
$$
donc
$$
\log N/\log \Im z \ioe \log \Im z \ioe \log \tau \sqrt{2}.
$$
D'autre part,
\begin{align*}
  V & \ioe  \demi  \log (n+1)/\log \log (n+1) +(1/2+\delta)\log (n+1) \log \log \log (n+1) /(\log \log (n+1))^2 +1\\
 & \ioe \demi  \log \tau/\log \log \tau +(1/2+2\delta)\log \tau \log \log \log \tau /(\log \log \tau)^2.
\end{align*}
Par conséquent,
\begin{align*}
V\log (\log N/\log \Im z) +(2+3\delta)V\log\log V  \ioe &\demi  (\log \tau/\log \log \tau) \log (\log N/\log \tau)\\
& + (3/2+5\delta)\log \tau \log \log \log \tau /\log \log \tau.
\end{align*}

On a donc montré que
$$
 \sup_{z \in \Scal_{N,\tau}}|\zeta(z)^{-1}N^{z}| \ioe \sqrt{N} \exp \Bigl (\demi  (\log \tau/\log \log \tau) \log (\log N/\log \tau) +(3/2+5\delta) \log \tau \log \log \log \tau /\log \log \tau \Bigr ).
$$
Ainsi, pour $\sqrt{2}N_0\ioe|\tau|\ioe\frac{1}{\sqrt{2}}N_1$, on a 
\begin{align*}
\intc_{\Scal_{N,\tau}}\zeta(z)^{-1}\frac{N^{z}}{z-i\tau}dz
\ioe &\sqrt{N} \exp \left( (\log |\tau|/2\log \log |\tau|) \log (\log N/\log |\tau|)\right. \\
& \left. + (3/2+6\delta)\log |\tau| \log \log \log |\tau| /\log \log |\tau| \right),
\end{align*} 
ce qui donne finalement, en utilisant \eqref{t46}
\begin{align}\label{t47}
B_N\ioe &\sqrt{N} \exp \Bigl (\demi  (\log |\tau|/\log \log |\tau|) \log (\log N/\log |\tau|) + (3/2+6\delta)\log |\tau| \log \log \log |\tau| /\log \log |\tau| \Bigr )\nonumber\\
&+O\left( \sqrt{N} \exp\bigl ((\log N)^{1/2}(\log \log N)^{5/2+6\delta}\bigr )\right).
\end{align}  

\subsubsection*{Conclusion : estimation de $M_N(i\tau)$}

D'après \eqref{t45} et \eqref{t47}, on a
$$
M_N(i\tau)=N^{-i\tau}B_N +O(|\tau|\log N) \quad (1 \ioe |\tau| \ioe N/5)
$$
et
\begin{align*}
B_N\ioe &\sqrt{N} \exp \Bigl (\demi  (\log |\tau|/\log \log |\tau|) \log (\log N/\log |\tau|) + (3/2+6\delta)\log |\tau| \log \log \log |\tau| /\log \log |\tau| \Bigr )\\
&+O\left( \sqrt{N} \exp\bigl ((\log N)^{1/2}(\log \log N)^{5/2+6\delta}\bigr )\right).
\end{align*}
On observe que sous les hypothèses de la proposition, on a :
$$|\tau|\log N \ioe N^{1/2}|\tau|^{2/5}$$
et
$$
N^{1/2}\exp\bigl ((\log N)^{1/2}(\log \log N)^{5/2+6\delta}\bigr ) \ioe N^{1/2}|\tau|^{1/3}.
$$
On a également
\begin{align*}
  \frac{\log |\tau|}{(\log \log |\tau|)^{5/2}} & \soe \frac{3(\log N)^{1/2}(\log \log N)^{5/2}}{\Bigl (\log  \bigl (3(\log N)^{1/2}(\log \log N)^{5/2}\bigr )\Bigr )^{5/2}}\\
& \soe \sqrt{\log N}. 
\end{align*}

Par conséquent,
$$
\frac{\log N}{\log |\tau|} \ioe \frac{\log |\tau|}{(\log \log |\tau|)^{5}},
$$
ce qui implique
\begin{align*}
\demi  \frac{\log |\tau|}{\log \log |\tau|} \cdot\log \Bigl (\frac{\log N}{\log |\tau|}\Bigr )  + (3/2+6\delta)\log |\tau|\frac{ \log \log \log |\tau|}{\log \log |\tau|} & \ioe \demi  \log |\tau|  +(-1+6\delta)\log |\tau| \frac{\log \log \log |\tau| }{\log \log |\tau|}
\end{align*}
et permet de conclure.\fin

\subsection{Estimations de $\zeta (s+\eps)^{-1}-M_N(s+\eps)$}

Démontrons à présent la proposition \ref{t52} et revenons à l'estimation de la différence
$$
\zeta (s+\eps)^{-1}-M_N(s+\eps),
$$
que nous exprimons d'abord à l'aide d'une intégrale :

  \begin{equation}
    \label{t49}
 \zeta (s+\eps)^{-1}-M_N(s+\eps)= -M_N(i\tau)N^{-1/2-\eps}+(1/2+\eps)\int_N^{\infty} t^{-3/2-\eps}M_t(i\tau)dt \quad (N \soe 1, \,\eps >0, \,\tau \in \Real)    
  \end{equation}
On suppose $N$ assez grand, $\eps \soe 2(\log \log N)^{5/2 +\delta}(\log N)^{-1/2}$, et $\tau \in \Real$.

\subsubsection*{Petites valeurs de $|\tau|$}

On a d'abord, d'après la proposition \ref{t48},
$$
M_N(i\tau)N^{-1/2-\eps} \ll (1+|\tau|)N^{-\eps}\exp\bigl ((\log N)^{1/2}(\log \log N)^{5/2 +\delta}\bigr ).
$$

D'autre part, pour $t\soe N$, on a
$$
\frac{\eps}{2}\log t \soe  (\log t)^{1/2}(\log \log t)^{5/2 +\delta}.
$$
En particulier,
$$
M_N(i\tau)N^{-1/2-\eps} \ll (1+|\tau|)N^{-\eps/2}.
$$

Et aussi,
\begin{align*}
\int_N^{\infty} t^{-3/2-\eps}M_t(i\tau)dt &\ll   (1+|\tau|)\int_N^{\infty} t^{-1-\eps}\exp\bigl ((\log t)^{1/2}(\log \log t)^{5/2 +\delta}\bigr )dt\\
&\ioe (1+|\tau|)\int_N^{\infty} t^{-1-\eps/2}dt\\
& \ll \eps^{-1} (1+|\tau|)N^{-\eps/2}.
\end{align*}

Or
\begin{align*}
  \eps^{-1} &\ioe (\log \log N)^{-5/2 }(\log N)^{1/2}\\
& \ioe \exp\bigl (\frac{1}{3}(\log N)^{1/2}(\log \log N)^{5/2 +\delta}\bigr )\\
& \ioe N^{\eps/6}, 
\end{align*}
donc $\eps^{-1} N^{-\eps/2} \ll N^{-\eps/3}$, ce qui donne sous nos hypothèses, la majoration
\begin{equation*}
\zeta (s+\eps)^{-1}-M_N(s+\eps) \ll (1+|\tau|)N^{-\eps/3}.
\end{equation*}
Dans le cas $\exp\bigl (3(\log N)^{1/2}(\log \log N)^{5/2 +6\delta}\bigr ) \soe |\tau|$, pour obtenir le résultat de la proposition \ref{t52}, il nous suffit donc de démontrer que
$$
 (1+|\tau|)N^{-\eps/3} \ll  (1+|\tau|)^{1/3}N^{-\eps/4},
$$
c'est-à-dire
$$
\frac{\eps}{12}\log N \soe \frac{2}{3}\log (1+|\tau|).
$$
Or on a bien dans ce cas
\begin{align*}
\frac{2}{3}\log (1+|\tau|) & \ioe \frac{2}{3} \bigl (3(\log N)^{1/2}(\log \log N)^{5/2 +6\delta}+O(1)\bigr )\\
&\ioe \frac{25}{12}(\log N)^{1/2}(\log \log N)^{5/2 +6\delta}\\
&\ioe \frac{\eps}{12}\log N. 
\end{align*}

\subsubsection*{Grandes valeurs de $|\tau|$}

Si $\exp\bigl (3(\log N)^{1/2}(\log \log N)^{5/2 +6\delta}\bigr ) \ioe |\tau| \ioe N^{3/4}$,
on a d'abord, d'après la proposition \ref{t44},
$$
M_N(i\tau)N^{-1/2-\eps} \ll N^{-\eps}|\tau|^{1/2-\kappa(\tau)}.
$$

\'Etudions maintenant l'intégrale
$$
\int_N^{\infty} t^{-3/2-\eps}M_t(i\tau)dt.
$$ 
Pour commencer, observons que $|\tau| \ioe N^{3/4} \ioe t^{3/4}$ si $t\soe N$.

D'autre part, définissons $\theta=\theta(\tau)$ par la relation
$$
|\tau| =\exp\bigl (3(\log \theta)^{1/2}(\log \log \theta)^{5/2 +6\delta}\bigr ).
$$
On a $\theta \soe N$ si $|\tau| \soe \exp\bigl (3(\log N)^{1/2}(\log \log N)^{5/2 +6\delta}\bigr )$, et
$$
\int_N^{\infty} t^{-3/2-\eps}M_t(i\tau)dt=\int_N^{\theta} t^{-3/2-\eps}M_t(i\tau)dt +\int_{\theta}^{\infty} t^{-3/2-\eps}M_t(i\tau)dt.
$$ 

Pour la première intégrale, nous pouvons utiliser la proposition \ref{t44} car $t \ioe \theta \implique |\tau| \soe \exp\bigl (3(\log t)^{1/2}(\log \log t)^{5/2 +6\delta}\bigr )$. Ainsi,
\begin{align*}
\int_N^{\theta} t^{-3/2-\eps}M_t(i\tau)dt & \ll |\tau|^{1/2-\kappa(\tau)}\int_N^{\theta} t^{-1-\eps}dt\\
& \ioe  |\tau|^{1/2-\kappa(\tau)}\eps^{-1}N^{-\eps}\\
& \ioe  |\tau|^{1/2-\kappa(\tau)}N^{-5\eps/6},\\  
\end{align*}
comme dans le cas précédent.

Pour la seconde intégrale, nous pouvons utiliser la proposition \ref{t48}. On a
$$
\int_{\theta}^{\infty} t^{-3/2-\eps}M_t(i\tau)dt \ll |\tau|\int_{\theta}^{\infty} t^{-1-\eps}  \exp\bigl ((\log t)^{1/2}(\log \log t)^{5/2 +6\delta}\bigr )dt.
$$

Maintenant, pour $t \soe \theta (\tau)$ $(\soe N)$, on a
$$
\frac{\eps}{2}\log t \soe 4(\log t)^{1/2}(\log \log t)^{5/2 +6\delta}.
$$
Ainsi,
\begin{align*}
\int_{\theta}^{\infty} t^{-3/2-\eps}M_t(i\tau)dt & \ll |\tau|\int_{\theta}^{\infty} t^{-1-\eps/2} \exp\bigl (-3(\log t)^{1/2}(\log \log t)^{5/2 +6\delta}\bigr ) dt\\
&\ioe |\tau|\exp\bigl (-3(\log \theta)^{1/2}(\log \log \theta)^{5/2 +6\delta}\bigr )\int_{\theta}^{\infty} t^{-1-\eps/2}dt\\
& =(2/\eps)\theta^{-\eps/2}\\
&\ioe (2/\eps)N^{-\eps/2}\\
& \ll N^{-\eps/3}  
\end{align*}
ce qui entraîne 
\begin{equation*}
\zeta (s+\eps)^{-1}-M_N(s+\eps) \ll N^{-\eps/3}|\tau|^{1/2-\kappa(\tau)}
\end{equation*}

Notons à présent que pour $|\tau|$ grand, on a $\beta(\tau)-\kappa(\tau) \ll 1/\log|\tau|$. Cela permet de conclure la démonstration de la proposition \ref{t52}.\fin

\section{Majoration de $I_{N,\eps}$}\label{t19}

Dans tout ce paragraphe, on pose $\sigma=\demi$, c'est-à-dire $s=\demi+i\tau$.

\begin{prop}(HR)
  Pour $N \soe 1$, $0< \eps \ioe 1/2$, on a
  \begin{equation}
    \label{t25}
\int_{
                |\tau| \soe N^{3/4}}
      |\zeta (s)|^2 |\zeta(s+\eps)-M_N(s+\eps)|^2\frac{d\tau}{|s|^2} \ll N^{-1/9}.
  \end{equation}
\end{prop}
\dem

Il suffit de démontrer que, pour $T \soe 1$,
\begin{equation}
  \label{t26}
I_N(T,\eps) =  \int_{
                T\ioe |\tau| \ioe 2T}
      |\zeta (s)|^2 |\zeta(s+\eps)^{-1}-M_N(s+\eps)|^2\frac{d\tau}{|s|^2} \ll T^{-3/2}(T+N)\log N,
\end{equation}
car \eqref{t25} résultera de la sommation de \eqref{t26} pour les valeurs $T=2^kN^{3/4}$, $k \in \Nat$.

On a
$$
 I_N(T,\eps) \ll T^{-2}\int_{
                T\ioe \tau \ioe 2T}
      |\zeta (s)/\zeta(s+\eps)|^2 d\tau +4T^{-2}\int_{
                T\ioe \tau \ioe 2T}
      |\zeta (s)|^2 |M_N(s+\eps)|^2 d\tau.
$$

D'une part,
$$
\int_{
                T\ioe \tau \ioe 2T}
      |\zeta (s)/\zeta(s+\eps)|^2 d\tau \ll T^{3/2},
$$
d'après le point \textit{(i)} de la proposition \ref{t1}.

D'autre part,
\begin{align*}
\int_{
                T\ioe \tau \ioe 2T}
      |\zeta (s)|^2 |M_N(s+\eps)|^2 d\tau  &\ioe T^{1/2}\int_{
                T\ioe \tau \ioe 2T} |\sum_{n \ioe N}\mu (n) n^{-1/2-\eps}n^{-i\tau}|^2 d\tau,
\end{align*}
d'après l'inégalité $|\zeta (s)| \ll \tau^{1/4}$ (cf. \cite{T1986}, (5.1.8) p.96).

La dernière intégrale vaut
$$
\bigl (T+O(N)\bigr ) \sum_{n \ioe N}\mu^2 (n) n^{-1-2\eps}\ioe (T+N)\log N,
$$
d'après une inégalité de Montgomery et Vaughan (cf. \cite {M1994}, (5) p.128), et car $ \sum_{n \ioe N}n^{-1-2\eps} \ll \log N$. Par conséquent,
\begin{equation*}
 I_{N,\eps} \ll  T^{-3/2} (T+N)\log N.\fine
\end{equation*}

\begin{prop}(HR)
Soit $N$ assez grand et $\eps \soe 25(\log \log N)^{5/2 +6\delta}(\log N)^{-1/2}$. Alors,
$$
\int_{|\tau| \ioe N^{3/4}} |\zeta (s)|^2 |\zeta(s+\eps)^{-1}-M_N(s+\eps)|^2\frac{d\tau}{|s|^2} \ll N^{-\eps/2}.
$$
\end{prop}
\dem  
 
Pour $|\tau| \ioe N^{3/4}$, on a
$$
\zeta (s+\eps)^{-1}-M_N(s+\eps) \ll N^{-\eps/4}(1+|\tau|)^{1/2-\beta(\tau)},
$$
d'après la proposition \ref{t52}. D'autre part,
\begin{align*}
  |\zeta (s)|^2 &\ll \exp\Bigl (O\bigl(\log (3+|\tau|)/\log \log (3+|\tau|)\bigr )\Bigr ) &\text{\footnotesize (\cite{T1986}, (14.14.1))}\\
& \ll (1+|\tau|)^{\beta(\tau)},
\end{align*}
donc
$$
\int_{|\tau| \ioe N^{3/4}} |\zeta (s)|^2 |\zeta(s+\eps)-M_N(s+\eps)|^2\frac{d\tau}{|s|^2} \ll
N^{-\eps/2}\int_{-\infty}^{\infty}(1+|\tau|)^{-1-\beta(\tau)}d\tau,
$$
où la dernière intégrale est convergente.\fin

\medskip

Les deux propositions précédentes entraînent la proposition \ref{t57}, ce qui achève la démonstration du théorème.

\medskip

\begin{multicols}{2}
\footnotesize  

\noindent BALAZARD, Michel\\
Institut de Math\'ematiques de Luminy, UMR 6206\\
CNRS, Université de la Méditerranée\\
Case 907\\
13288 Marseille Cedex 09\\
FRANCE\\
Adresse \'electronique : \texttt{balazard@iml.univ-mrs.fr}

\smallskip

\noindent de ROTON, Anne\\
Institut Elie Cartan de Nancy, UMR 7502\\
Nancy-Université, CNRS, INRIA\\
BP 239\\
54506 Vandoeuvre-lès-Nancy Cedex\\
FRANCE\\
Adresse \'electronique : \texttt{deroton@iecn.u-nancy.fr
}
\end{multicols}

\end{document}